\documentclass[12pt]{amsart}
\usepackage{amscd,amssymb,subfigure}
\usepackage[arrow,matrix,graph,frame,poly,arc]{xy}
\usepackage{pstricks}
\usepackage{graphicx}

\newtheorem{theorem}[subsection]{Theorem}
\newtheorem{lemma}[subsection]{Lemma}
\newtheorem{prop}[subsection]{Proposition}
\newtheorem{corollary}[subsection]{Corollary}

\newtheorem{thm}{Theorem}

\theoremstyle{definition}
\newtheorem{remark}[subsection]{Remark}
\newtheorem{definition}[subsection]{Definition}
\newtheorem{example}[subsection]{Example}

\numberwithin{equation}{section}
\newcommand{\A}{{\mathcal A}}

\renewcommand{\L}{{\mathcal L}}

\newcommand{\Z}{\mathbb{Z}}
\newcommand{\Q}{\mathbb{Q}}

\newcommand{\C}{\mathbb{C}}
\newcommand{\F}{\mathbb{F}}
\newcommand{\TT}{\mathbb{T}}

\renewcommand{\k}{\Bbbk}

\newcommand{\G}{\Gamma}
\newcommand{\E}{\mathsf{E}}
\newcommand{\V}{\mathsf{V}}

\DeclareMathOperator{\Hom}{Hom}
\DeclareMathOperator{\rk}{rk}

\DeclareMathOperator{\ch}{char}

\newcommand{\abs}[1]{\left| #1 \right|}

\begin{document}

\title[Milnor fibers of graphic arrangements]{%
On the monodromy action on Milnor fibers of graphic arrangements}

\author[A.~D.~Macinic]{Anca Daniela M\u acinic$^*$}
\address{Inst. of Math. Simion Stoilow,
P.O. Box 1-764,
RO-014700 Bucharest, Romania}
\email{Anca.Macinic@imar.ro}

\author[S.~Papadima]{\c Stefan~Papadima$^*$}
\address{Inst. of Math. Simion Stoilow,
P.O. Box 1-764,
RO-014700 Bucharest, Romania}
\email{Stefan.Papadima@imar.ro}

\thanks{$^*$Partially supported by the CEEX Programme of
the Romanian Ministry of Education and Research, contract
2-CEx 06-11-20/2006.}



\begin{abstract}
We analyze the monodromy action, over the rationals, on the first
homology group of the Milnor fiber, for arbitrary subarrangements
of Coxeter arrangements. We propose a combinatorial formula for
the monodromy action, involving Aomoto complexes in positive
characteristic. We verify the formula, in cases $A$, $B$ and $D$.
\end{abstract}

\maketitle


\section{Introduction and statement of results}
\label{sect:intro}

Let $\A =\{ H_1,\dots, H_n \}$ be an arrangement of complex hyperplanes in
 $\C^{l}$,
with {\em complement} $M_{\A}= \C^{l} \setminus \bigcup_{i=1}^n H_i$, 
and {\em intersection lattice} $\L (\A)$, consisting of the various 
intersections of hyperplanes from $\A$, ordered by reverse inclusion. A
fundamental result in arrangement theory, due to Orlik and Solomon \cite{OS},
relates the topology and the combinatorics of $\A$, by saying that the homology
of $M_{\A}$ with arbitrary untwisted coefficients is {\em combinatorial},
i.e., is determined by the intersection lattice. More precisely, they proved
that the cohomology ring with arbitrary coefficients, $H^*(M_\A, \k)$,
is isomorphic to the Orlik-Solomon algebra of $\A$ over $\k$, $A_{\k}^*(\A)$,
which depends only on the lattice $\L(\A)$. 

Assuming $\A$ to be central, with homogeneous defining polynomial, $f_{\A}$,
there is a well-known global Milnor fibration, $F_{\A}\hookrightarrow M_{\A}
\stackrel{f_{\A}}{\rightarrow}\C^*$, where $F_{\A}:=f_{\A}^{-1}(1)$ is the
{\em Milnor fiber}. Milnor fibers of polynomials and their homology, especially
the structure of the {\em monodromy action} on homology, play a key role in
singularity theory, see for instance \cite{D} and the references therein.
An important problem in arrangement theory is to
decide whether $H_*(F_{\A}, \Q)$ is combinatorially determined. 
To our best knowledge, the problem is
open, even in degree $*=1$. (Libgober's description \cite{Li1, Li2} of the
monodromy action, in terms of superabundance of curves, is apriori non-combinatorial.)

The finite {\em graphs} $\G$ we consider in this paper, with vertex set $\V$
and edges $\E$, have at most double edges connecting two distinct vertices, and
at most one loop at each point. The presence of a loop at $i$ will be
denoted by $\bigodot i$. Edges are labeled with signs: double edges are 
indicated by the label $\pm$, positive simple edges by $+$, and the absence of
a label indicates a negative edge. 

An {\em unsigned graph} means an ordinary
finite simplicial graph (with no double edges or loops), where all edges are
negative. A {\em signed graph} is a graph without loops. The graphs $\G$ we are
considering here encode subarrangements of Coxeter arrangements of type $B$, called
{\em graphic arrangements} and denoted by $\A(\G)$. The signed graphs correspond to
subarrangements of  Coxeter arrangements of type $D$, while the unsigned ones 
parametrize type $A$ subarrangements. The definition of $\A(\G)$ is the obvious one;
see Definition \ref{def:agraph}.

For example, in the figure below $\G$ is unsigned, whereas $\G'$ has a double edge,
$5$ negative edges, $4$ positive edges, and $3$ loops.

\begin{pspicture}(0,1.5)
\pscircle*(3,1){.07}
\pscircle*(5,1){.07}
\pscircle*(4,0){.07}
\pscircle*(3,-1){.07}
\pscircle*(5,-1){.07}
\psline[linewidth=.5pt](5,1)(5,-1)(3,-1)(3,1)(5,1)(4,0)(3,-1)
\psline[linewidth=.5pt](3,1)(4,0)(5,-1)
\pscircle*(8,-1){.07}
\pscircle*(8,1){.07}
\pscircle*(10,1){.07}
\pscircle(10,1){.4}
\pscircle*(9,0){.07}
\pscircle*(10,-1){.07}
\pscircle(10,-1){.4}
\pscircle*(11,0){.07}
\pscircle(11,0){.4}
\psline[linewidth=.5pt](10,1)(10,-1)(8,-1)(8,1)(10,1)(9,0)(10,-1)
\psline[linewidth=.5pt](8,1)(9,0)(8,-1)

\rput(9.6,.4){\footnotesize $+$}
\rput(8.2,0){\footnotesize $+$}
\rput(10.2,0){\footnotesize $+$}
\rput(9,1.2){\footnotesize $\pm$}
\rput(2.7,0){\footnotesize $\Gamma$}
\rput(7.6,0){\footnotesize $\Gamma'$}
\end{pspicture}
\vskip .5in

\begin{figure}[ht]
\caption{\textsf{Two graphs}}
\label{fig:graph}
\end{figure}

Since the geometric monodromy action on $F_{\A}$ has order $n$, it follows that
one has an equivariant decomposition (with respect to the homology monodromy
action),
\begin{equation}
\label{eq:ciclo}
H_q(F_{\A}, \Q)= \bigoplus_{d |n} \big( \frac{\Q[t]}{\Phi_d} \big)^{b_{qd}(\A)}
\end{equation}
for all $q$, where $\Phi_d$ is the $d$th cyclotomic polynomial; see
\cite{OT, L}.

The numbers $b_{q1}(\A)$, $q\ge 0$, are combinatorially determined, being
equal to the
corresponding Betti numbers of the associated projective arrangement 
$\overline{\A}$;
see \cite{OT}. In particular, $b_{11}(\A)= n-1$. We may also assume in 
\eqref{eq:ciclo} that $r:= \rk (\A)\ge 3$ (if $r=1$, $F_{\A}$ is a point, 
and the rank
$2$ case is treated in \cite[Proposition 5.125]{OT}). 

Our main result in this paper establishes a combinatorial formula for
the numbers
$b_d(\G):= b_{1d}(\A(\G))$, in the case of graphic arrangements. To
describe it,
we need to recall the general definition of {\em Aomoto complexes} associated
to Orlik-Solomon algebras, $A_{\k}^*(\A)$. Let $\omega \in A_{\k}^1(\A)$ be an
arbitrary
element. Since $A^*$ is a quotient of an exterior algebra, the square 
$\omega \cdot \omega$ vanishes. Denoting by $\mu_{\omega}$ left-multiplication
by $\omega$ in $A^*$, we thus obtain a cochain complex,
\begin{equation}
\label{eq:aok}
\big( A_{\k}^*(\A), \mu_{\omega} \big)= \{ A_{\k}^*(\A)\stackrel{\mu_{\omega}
}{\longrightarrow} 
A_{\k}^{*+1}(\A)\}_{*\ge 0}\, ,
\end{equation}
called the Aomoto complex of $\omega$, introduced by Aomoto in \cite{A},
and studied by Falk in \cite{F}, from the point of view of 
resonance varieties of arrangements.

By definition, $A_{\k}^1(\A)$ is freely generated by $\{ a_H \}_{H\in \A}$. So,
$\omega= \sum_{H\in \A} \lambda_H a_H$, with $\lambda_H\in \k$.
Denote by $\omega_1$ the
distinguished element $\omega_1:= \sum_{H} a_H$, and abbreviate
 $\mu_{\omega_1}$ by
$\mu_1$. Let $\k$ be  a field, $\ch \k =p$. Set
\begin{equation}
\label{eq:modres}
\beta_{qp}(\A):= \dim_{\k} H^q (A_{\k}^*(\A), \mu_1)\, \quad {\rm for}\quad 
q\ge0 \,.
\end{equation}
One knows \cite{Y} that $\beta_{q0}(\A)=0$, for all $q$. When $\A=\A(\G)$ is
a graphic arrangement, set $\beta_p(\G):= \beta_{1p}(\A(\G))$, for each prime $p$. 

The off-diagonal elements different from $2$ of type $A-I$ Coxeter matrices,
in rank $\ge 3$, are 
$3$, $4$ and $5$ \cite{B}. All of them are of the form $p^k$, with
 $p\in \{ 2,3,5 \}$.
The theorem below relates the numbers $b_d(\G)$ from \eqref{eq:ciclo} to the
numbers $\beta_p(\G)$ coming from \eqref{eq:aok}.

\begin{thm}
\label{thm:main}
Let $\A(\G)$ be an arbitrary graphic arrangement of rank at least $3$, with $n$
hyperplanes, and let $d\ne 1$ be a divisor of $n$. 
\begin{enumerate}
\item \label{a1} 
If $d\ne 3$, then $b_d(\G)=0$.
\item \label{a2}
If $p\ne 3$ is prime, then $\beta_p(\G)=0$.
\item \label{a3}
If $n\equiv 0$ (mod $3$), then $b_3(\G)=\beta_3(\G)$. 
If $n \not\equiv 0$ (mod $3$), then $\beta_3(\G)=0$.
\item \label{a4}
The following formula holds for the Milnor fiber $F_{\G}$:
\[
H_1(F_{\G}, \Q)= \big( \frac{\Q[t]}{t-1} \big)^{n-1}\oplus
\big( \frac{\Q[t]}{\Phi_2} \oplus
\frac{\Q[t]}{\Phi_4} \big)^{\beta_{2}(\G)} \oplus
\big( \frac{\Q[t]}{\Phi_3} \big)^{\beta_{3}(\G)}\oplus
\big( \frac{\Q[t]}{\Phi_5} \big)^{\beta_{5}(\G)}\, .
\]
\end{enumerate}
\end{thm}

We conjecture that the above formula \eqref{a4} actually holds for all
subarrangements $\A$ of rank at least $3$ of arbitrary Coxeter arrangements.

A useful fact is that the $\Q [t]$--module structure of $H_*(F_{\A}, \Q)$ 
depends only on the lattice-isotopy type (in the sense of Randell \cite{R}) of
the arrangement $\A$; see Section \ref{sect:cyclolocal} for more details. 
With this remark,
\eqref{eq:ciclo} takes the following explicit form, when $\A$ is graphic.

\begin{thm}
\label{thm:maincompl}
Let $\A(\G)$ be an arbitrary graphic arrangement of rank at least $3$, with
Milnor fiber $F_{\G}$. Set $n:= \abs{\E(\G)} $.
\begin{enumerate}
\item \label{b1}
If $\A(\G)$ is lattice-isotopic to either $D_3$ or $D_4$ (the full Coxeter
arrangements of type $D$ and rank $3$ or $4$), then
\[
H_{1}(F_{\G}, \mathbb{Q})=\big( \frac{\mathbb{Q}[t]}{t-1} \big)^{n-1}\oplus
\big( \frac{\mathbb{Q}[t]}{t^{2}+t+1} \big)\, .
\]
\item \label{b2}
Otherwise, $H_{1}(F_{\G}, \mathbb{Q})=\big( \frac{\mathbb{Q}[t]}{t-1} \big)^
{n-1}$.
\end{enumerate}
\end{thm}

Similar results (proving the asymptotic triviality of the monodromy action on
$H_{q}(F_{\G}, \mathbb{Q})$) have been obtained by Settepanella, in the 
particular case of complete graphic arrangements on $v\gg q$ vertices, of types
$A$, $B$ and $D$; see \cite{S}.

However, the methods are completely different. The main tool from \cite{S} is
the Salvetti complex associated to a Coxeter group. This technique does not seem
to extend to arbitrary subarrangements of Coxeter arrangements. Our strategy is
to use the known relationship between Milnor fibers and twisted homology, see
for instance Cohen-Suciu \cite{CS}. To compute the latter, via Aomoto complexes,
we rely on three key results: the first in characteristic zero 
(\cite{ESV, STV}),
the second in arbitrary characteristic (\cite{Y}), and the last in 
positive characteristic (\cite{OC, PS}). These techniques are available for 
arbitrary arrangements $\A$.

Based on a method due to Deligne \cite{De}, Esnault-Schechtman-Viehweg \cite{ESV}
and Schechtman-Terao-Varchenko \cite{STV} showed that twisted homology on $M_{\A}$
may be computed by Aomoto complexes in characteristic zero, 
for certain local systems. Unfortunately, this
approach does not always work, see e.g. Example \ref{ex=k6}. When the Deligne
method is available, it may be combined with general results on Aomoto 
complexes, due
to Yuzvinsky \cite{Y}, to obtain vanishing results. We use this approach in
Theorem \ref{thm:main} \eqref{a1}, for $d\ne 2,3,4$. 

To settle the remaining cases, we resort to {\em modular upper bounds}, for the
dimension over $\C$ of twisted homology with rational local systems whose
denominator is a prime power, $p^k$. Improving a result due to Cohen and Orlik 
\cite{OC} for $k=1$, 
it is shown in \cite{PS} that these dimensions are bounded above by numbers
coming from objects in characteristic $p$; in the equimonodromical case, these
numbers are defined by \eqref{eq:modres}. This method yields 
Theorem \ref{thm:maincompl} \eqref{b2}.

In all previously known (sporadic) examples, the modular inequalities become
equalities, for equimonodromical rational local systems with $k=1$; 
see \cite[Section 7]{CO1}. We
may add the following new large class of examples to the list.

\begin{thm}
\label{thm:sharp}
Let $\A(\G)$ be a graphic arrangement (of arbitrary rank). The modular upper 
bound for equimonodromical rational local systems on $M_{\A(\G)}$ 
with denominator $p$ is equal to
the dimension of the corresponding twisted homology in degree one, for every
prime $p$.
\end{thm}

Our approach also leads to a partial verification of formula \eqref{a4} from
Theorem \ref{thm:main}, for arbitrary subarrangements of arbitrary Coxeter type;
see Corollary \ref{cor:conj}.

\section{Homology of Milnor fibers and twisted homology}
\label{sect:cyclolocal}

In this section, we will review the relationship between the
cyclotomic decomposition of $H_*(F_{\A}, \Q)$, and the
(co)homology of the complement of $\A$ with coefficients in rank one 
local systems.

Assume $\mathcal{A}$ is an arrangement in $\mathbb{C}^{l}$, defined as 
the zero set of the homogeneous polynomial
$f_{\mathcal{A}}$.
There is an action on $\mathbb{C}^{l}$, given by the multiplication with 
$u=\exp\frac{2\pi \sqrt{-1}}{n}$, where $n=|\mathcal{A}|$,  which induces 
an action on the fiber $F_{\mathcal{A}}$ 
(since $f_{\mathcal{A}}$ is homogeneous of degree $n$). 
We call this action 
on the Milnor fiber the \textit{geometric monodromy}, denoted by 
$h:F_{\mathcal{A}}\longrightarrow F_{\mathcal{A}}$. The induced action on
homology, $h_* \colon H_*(F_{\A}, \Q)\to H_*(F_{\A}, \Q)$, corresponds to
multiplication by $t$, in equation \eqref{eq:ciclo}.

\subsection{}
\label{ss21}

This may be conveniently reinterpreted in terms of {\em twisted homology},
as follows. The complement $M_{\A}$ is a connected, $1$--{\em marked}, 
finite type CW--space.
That is, it is endowed with a $\Z$--basis of $H_1(M_{\A})$, denoted by
$\{ a^*_H \}_{H\in \A}$, dual to the canonical basis of $A^1_{\Z}(\A)$. The
marking defines a $\Z$--character, $\nu \colon H_1(M_{\A})\to \Z$, which sends
each $a^*_H$ to $1$. This character induces on group rings a homomorphism,
$\nu \colon \Z \pi_1(M_{\A})\to \Z[t^{\pm 1}]$, which gives rise to a
$\Z \pi_1(M_{\A})$--module (alias a local system on  $M_{\A}$), denoted by
$\Q[t^{\pm 1}]_{\nu}$. There is an equivariant isomorphism
\begin{equation}
\label{eq:cyclicf}
H_*(F_{\A}, \Q)\cong H_*(M_{\A}, \Q[t^{\pm 1}]_{\nu})\, ,
\end{equation}
see \cite[p.106--107]{D} and \cite[Ch.VI]{W}.

\subsection{}
\label{ss22}

One may consider arbitrary ring homomorphisms
$\nu \colon \Z \pi_1(M_{\A})\to R$, where $R$ is a commutative ring,
with group of units $R^*$. These morphisms are naturally identified with elements of
$\Hom (H_1(M_{\A}), R^*)\equiv (R^*)^n$. The associated local system, $R_{\nu}$,
is called {\em equimonodromical} if $\nu$ is constant on the distinguished basis
$\{ a^*_H \}$. It follows from \cite[p.497-498]{P} that the equivariant isomorphism
type of $H_*(M_{\A}, R_{\nu})$ depends only on the lattice-isotopy type of $\A$, in
the equimonodromical case. From the definitions, we also see that
the cochain isomorphism type of the Aomoto complex
$(A_{\k}^*(\A), \mu_1)$ defined in the Introduction depends only on lattice-isotopy type.

\subsection{}
\label{ss23} 

Twisted homology with coefficients in rank one local systems,
$H_*(M_{\A}, \C_{\rho})$, is a very active research area in arrangement theory.
Here, $\rho\in \Hom (H_1(M_{\A}), \C^*)$ denotes an arbitrary {\em character}.
The {\em rational} characters play an important role.

\begin{definition}
\label{def:ratl}
Let ${\bf k} =(k_H)_{H\in \A}$ be a collection of integers, with g.c.d. equal to $1$.
Let $u\in \C^*$ be a primitive $d$--root of unity. The character $\rho$ defined by
$\rho (a^*_H)= u^{k_H}$ is called rational. If $\bf{k}= \bf{1}$, $\rho$ is called 
rational and equimonodromical, with denominator $d$.
\end{definition}

Set $b_q(\A, \frac{\bf{k}}{d}):= \dim_{\C} H_q(M_{\A}, \C_{\rho})$. (This is
well-defined, by Galois theory.) As is well-known (see e.g. \cite{DN}), one has
the following recurrence formula, for $d\mid n$:

\begin{equation}
\label{eq:twistedf}
b_q(\A, \frac{\bf{1}}{d})= b_{qd}(\A) +b_{q-1,d}(\A)\, , \forall q \,.
\end{equation}
In particular, $b_{d}(\A):= b_{1d}(\A)= b_1(\A, \frac{\bf{1}}{d})$, for
$1\ne d \mid n$. 

\subsection{}
\label{ss24}

We close this section by describing a method for computing twisted homology on $M_{\A}$,
by using {\em generic sections}.

We will need the following version of twisted Betti numbers, for arbitrary Aomoto complexes.
Given $\omega \in A^1_{\k}(\A)$, $\k$ a field, set

\begin{equation}
\label{eq:baom}
\beta_q(\A, \omega):= \dim_{\k} H^q (A^{\bullet}_{\k}(\A), \mu_{\omega}) 
\quad {\rm for} \quad q\ge 0 \,.
\end{equation}
We may now spell out our result.

\begin{prop}
\label{prop:gens}
Let $\A$ be a rank $r\ge 3$ arrangement in $\C^l$. Let $U\subset \C^l$ be a subspace of
dimension $k+1$, $2\le k<r$. Denote by $\A^U$ the restriction, and by
$j\colon M_{\A}\cap U\hookrightarrow M_{\A}$ the inclusion map between complements.
If $U$ is $\L_k(\A)$--generic, in the sense of \cite[\S 5(1)]{PD}, the following hold.
\begin{enumerate}
\item \label{g0}
The map induced by $j$ on $\pi_1$ is an isomorphism, preserving the natural $1$--markings
upon abelianization.

\item \label{g1}
The map induced on Aomoto complexes,
$j^*\colon (A_{\k}^*(\A), \mu_{\omega})\to (A_{\k}^*(\A^U), \mu_{\omega})$,
is an isomorphism for $*\le k$. In particular, 
$\beta_q(\A, \omega)= \beta_q(\A^U, \omega)$, for any  $\omega$ and every $q<k$.

\item \label{g2}
The map induced on twisted homology,
$j_*\colon H_*(M_{\A}\cap U, j^*R)\to H_*(M_{\A}, R)$, is an isomorphism for
$*<k$ and an epimorphism for $*=k$, for arbitrary coefficients. 
Moreover, $j^*R\equiv R$, if $R$ comes from
a representation, $\nu \colon \Z \pi_1(M_{\A})\to R$, where $R$ is a commutative ring.

\end{enumerate}
\end{prop}

\begin{proof}
By \cite[Proposition 5.14]{PD}, $j$ induces an isomorphism on $\pi_q$, for $q<k$,
and a surjection on $\pi_k$. 

\eqref{g0} Remember that $k\ge 2$, to obtain the assertion on $\pi_1$. The claim on
markings is obvious. Put together, these two properties show that 
$j^*R\equiv R$, if $R$ comes from
an abelian representation.

\eqref{g1} Follows from the fact that $\A^U$ and $\A$ have the same 
dependent subarrangements,
up to cardinality $k+1$.

\eqref{g2} The first claim is a standard consequence of the properties of $j_{\sharp}$ on
$\pi_{\le k}$, see \cite[Ch.VI]{W}, and the second was already clarified in the proof of
\eqref{g0}.
\end{proof}

We will prove that, for almost all graphic arrangements,
the only nontrivial component  from decomposition \eqref{eq:ciclo} in degree $1$ 
is the part corresponding to $\Phi_1$.  
To do this, we turn to combinatorial computations.

\section{Twisted homology and Aomoto complexes}
\label{sect:generic}

Let $\omega \in A_{\C}^1(\A)$ be a degree one element of the Orlik--Solomon algebra
of $\A$ with complex coefficients. Write $\omega =\sum_{H\in \A}\lambda_H a_H$, with
$\lambda_H\in \C$. Consider the character torus, $\TT_{\A}:= \Hom (\pi_1(M_{\A}), \C^*)=
\Hom (H_1(M_{\A}), \C^*)\equiv (\C^*)^n$, and the rank one complex local system associated to 
$\omega$, $\rho_{\omega}:=(\exp (2\pi \sqrt{-1}\lambda_H))_{H\in \A}\in \TT_{\A}$.
Clearly, $\rho_{\omega}=\rho_{\omega +\alpha}$, for all $\alpha \in \Z^n$.

\subsection{}
\label{ss31}

Basic results from \cite{ESV, STV} establish a deep connection between the twisted 
cohomology of $M_{\A}$, $H^*(M_{\A}, {}_{\rho_{\omega}} \C)$, and the cohomology of the 
Aomoto complex of $\omega$,
$(A^{\bullet}_{\C}(\A), \mu_{\omega})$,
for {\em nonresonant} $\omega$.

\begin{definition}
\label{def=dense}
An element $X\in \L(\A)$ is called \textit{dense} if the arrangement
$\mathcal{A}_{X}$ is not
decomposable  as a nontrivial product.
\end{definition}

\begin{example}
\label{ex:dense12}

\textit{(i)} All hyperplanes are dense elements.

\noindent \textit{(ii)} An element $X$ of rank $2$  is dense if and only if
$|\mathcal{A}_{X}|\geq3$.
\end{example}

For $X\in \L(\A)$, set $m_X:= \abs{\A_X}$. For 
$\omega =\sum_{H\in \A}\lambda_H a_H \in A^1_{\k}(\A)$ and $X\in \L(\A)$, set
$\omega_X:= \sum_{H\supset X}\lambda_H a_H \in A^1_{\k}(\A_X)$, and
$\Sigma_X \omega:= \sum_{H\supset X}\lambda_H \in \k$. For a central arrangement $\A$,
let $C:= \cap_{H\in \A} H$ be the center of $\A$. 
 
\begin{definition}
\label{def=gen}
Let $\A$ be a central arrangement.
An element $\omega=\sum_{H\in{\mathcal{A}}}\lambda_{H}a_{H}\in{A^{1}_{\C}
(\mathcal{A})}$ is called \textit{nonresonant} if  
$\Sigma_X \omega\notin{\mathbb{Z}_{>0}}$,
for all dense elements 
$X\in{\L (\mathcal{A})}$,
and $\Sigma_C \omega\notin{\mathbb{Z}_{<0}}$.
\end{definition}

One may reduce the computation of twisted homology to a combinatorial problem, 
under a nonresonance assumption, via the following result.

\begin{theorem}[\cite{ESV, STV}]
\label{thm:gen}
Let $\omega\in{A^1_{\mathbb{C}}(\mathcal{A})}$ be a nonresonant element. Then 
$$\dim_{\C} H_{q}(M_{\mathcal{A}}, \mathbb{C}_{\rho_{\omega}})=
\beta_q(\A, \omega)\, , \forall q\, . $$
\end{theorem}

\subsection{}
\label{ss32}

We define now a partial nonresonance condition.

\begin{definition}
\label{def=kgen}
An element
$\omega=\sum_{H\in{\mathcal{A}}}\lambda_{H}a_{H}\in{A^1_{\C}(\mathcal{A})}$ is called
$k$--\textit{nonresonant} $(k\geq{1})$, if $\Sigma_X \omega \notin{\mathbb{Z}_{>0}}$,
for all dense elements $X\in{\L (\mathcal{A})}$ of rank $\leq{k+1}$, and
$\Sigma_C \omega =0$.
\end{definition}

This definition  leads to a refinement of Theorem \ref{thm:gen}.
 
\begin{prop}
\label{thm:partial gen}
Let $\A$ be a central arrangement, of rank $r\ge 3$.
If  $\omega \in A^{1}_{\C}(\mathcal{A})$ is $k$--nonresonant,
$1\le k<r-1$, then
\begin{equation}
\label{eq:knonres}
\dim_{\C} H_{q}(M_{\mathcal{A}}, \mathbb{C}_{\rho_{\omega}})=
\beta_q(\A, \omega)\, , \forall \, q\le k\, . 
\end{equation}
\end{prop}

\begin{proof}
Pick a subspace $U$, $(k+2)$--dimensional and $\L_{k+1}(\A)$--generic.
By Proposition \ref{prop:gens}, we may replace $\A$ by $\A^U$ in
\eqref{eq:knonres} above. Note also that $\rk (\A^U)=k+2$. Once we have
checked that $\omega \in A^{1}_{\C}(\mathcal{A}^U)$ is nonresonant, our claim
follows from Theorem \ref{thm:gen}.

To do this, we start by observing that the correspondence 
$X \leadsto X\cap U$ gives a bijection between $\L(\A)$ and $\L(\A^U)$,
in rank $\le k+1$. This is a direct consequence of the fact that $U$ is
$\L_{k+1}(\A)$--generic. Moreover, it is straightforward to verify that
this bijection is order and rank preserving, and induces a bijection 
$\A_X \stackrel{\sim}{\to} \A_{X\cap U}^U$, if $\rk (X)\le k+1$.

To check that the bijection also preserves dense elements, it is enough to recall
from \cite[Theorem 2]{C} that $X\in \L(\A)$ is dense if and only if
$(1+t)^2$ does not divide the Poincar\' e polynomial $P_{\A_X}(t)$.

Finally, just note that the partial nonresonance conditions for $\A$
coincide with the nonresonance conditions for $\A^U$, in rank $\le k+1$,
while the remaining nonresonance condition(s), for the center of $\A^U$,
take(s) a stronger form in $\A$; compare Definitions \ref{def=gen} and
\ref{def=kgen}. 
\end{proof}

\subsection{}
\label{ss33}

We would like to apply the above proposition to
$\frac{\bf{1}}{d}:= \sum_{H\in \A} \frac{a_H}{d}$. But the $1$--nonresonance
condition is clearly violated, as soon as $X$ has rank $2$, $m_X>2$ 
and $d \mid m_X$; see
Example \ref{ex:dense12}$(ii)$. This prompts the next definition.

\begin{definition}
\label{def:kadm}
An element $\omega \in A^{1}_{\C}(\mathcal{A})$ is $k$--{\em admissible} if there is
$\alpha \in \Z^n$ such that $\omega +\alpha$ is $k$--nonresonant.
\end{definition}

\begin{corollary}
\label{cor:kadm}
Assume $\rk (\A)\ge 3$. Let $\rho\in \TT_{\A}$ be a rational character. If
$\frac{\bf{k}}{d}$ is $k$--admissible, then $b_q(\A, \frac{\bf{k}}{d})=
\beta_q(\A, \frac{\bf{k}}{d}+\alpha)$, $\forall\, q\le k$, where $\alpha$ is 
as in Definition \ref{def:kadm}.
\end{corollary}

Unfortunately, there are simple nonadmissible examples.

\begin{example}
\label{ex=k6}
Let $\A$ be the full Coxeter arrangement $A_{v-1}$, corresponding to the
complete unsigned graph on $v$ vertices. When $v\ge 5$, 
$\frac{\bf{1}}{3}$ is not $1$--admissible.

Assuming the contrary, we may find $\alpha_{ij}\in \Z$, $1\le i\ne j\le v$, 
with the property that $\alpha_{ij}+ \alpha_{jk}+ \alpha_{ki}\ge 1$, for all
distinct $i, j, k$, and $\Sigma_C \alpha =\frac{v(v-1)}{6}$. Summing the 
inequalities, we get $(v-2)\Sigma_C \alpha\ge \binom{v}{3}$. Therefore, all
inequalities must be equalities. Solving the system for $v=4$, we find out that
necessarily $\alpha_{ij}=\alpha_{kl}$, if $i,j,k,l$ are distinct. If $v\ge 5$,
this implies that $\alpha_{ij}= \alpha_{jk} =\alpha_{ki} =\frac{1}{3}$, for all
distinct $i, j, k$, a contradiction. 
\end{example}

Nevertheless, Corollary \ref{cor:kadm} turns out to be very useful to obtain vanishing
results. To make this statement precise, we need the following definitions. For a given
arrangement $\A$ and for each $k\ge 2$, set
\begin{equation}
\label{eq:mlist}
{\bf m}_k(\A):= \{ m_X\, \mid \, X\in \L(\A) \quad {\rm dense}\quad {\rm and}\quad
2\le \rk(X)\le k \} \,.
\end{equation}
For a fixed hyperplane $K\in \A$, set also
\begin{equation}
\label{eq:multawayk}
{\bf m}_k^K (\A):= \{ m_X\, \mid \, X\in \L(\A) \quad {\rm dense}\, ,\,\,\, X\not\subset K
\quad 
{\rm and}\quad
2\le \rk(X)\le k \} \,.
\end{equation}
We may now state our result.

\begin{theorem}
\label{thm:negvanish}
Let $\A$ be a central arrangement of rank $r\ge 3$, with $n$ hyperplanes, and $1\le k<r-1$.
If $1\ne d\mid n$ is such that $d$ does not divide $m$, for any $m\in {\bf m}_{k+1}^K (\A)$,
for some $K\in \A$, then $b_{qd}(\A)=0$, for all $q\le k$.
\end{theorem}

\begin{proof}
Define $\alpha \in \Z^n$ by: $\alpha_H= 0$ (for $H\ne K$),
and $\alpha_K= \frac{n}{d}$. We claim that $\omega:= \frac{{\bf 1}}{d} -\alpha $ 
is $k$--nonresonant.
Plainly, $\Sigma_C \omega= 0$. The rank one nonresonance conditions involve
$\Sigma_H \omega$, which equals either $\frac{1}{d}$ (if $H\ne K$), which is not an integer, 
or $\frac{1-n}{d} <0$ (if $H=K$). For $X$ dense, $X\not\subset K$, with $2\le \rk(X)\le k+1$,
$\Sigma_X \omega= \frac{m_X}{d}- \Sigma_X \alpha$ cannot be an integer, 
since $d$ does not divide $m_X$. If $X\subset K$, then $\Sigma_X \omega= (m_X-n)/d \le 0$.
Thus, the $k$--nonresonance claim is established.

Hence, Proposition \ref{thm:partial gen} applies, and guarantees that 
$b_q(\A, {\bf 1}/d)=
\beta_q(\A, \omega)$, for all $q\le k$. Our next claim is that 
$\beta_q(\A, \omega)=0$, if $q\le k$. This may be seen by using
\cite[Theorem 4.1(ii)]{Y}, as follows. Pick a $(k+2)$--subspace $U$, which is
$\L_{k+1}(\A)$--generic. Due to Proposition \ref{prop:gens} \eqref{g1}, we may replace
$\A$ by $\A^U$. 

Let us check now, for $\A^U$, the hypotheses needed in 
the abovementioned result of Yuzvinsky. As we have seen before, $\Sigma_C \omega=0$.
The remaining conditions involve $\Sigma_X \omega$, for $X\in \L(\A^U)$ with
$1\le \rk(X)\le k+1$. Recall from the proof of Proposition \ref{thm:partial gen} that
these elements $X$ are identified with the  elements $X$ from $\L(\A)$ 
of rank at most $k+1$; moreover, $\Sigma_X \omega$ takes the same value in $\A^U$ as in $\A$.

There are two cases to be considered. If $X\subset K$, then $\Sigma_X \omega= (m_X-n)/d<0$
(since $\A_{X\cap U}^U \ne \A^U$). Otherwise, $\Sigma_X \omega= m_X/d>0$. In both cases,
$\Sigma_X \omega\ne 0$, and we are done.

We may conclude by deducing inductively from $b_q(\A, {\bf 1}/d)=0$, for $q\le k$, that
$b_{qd}(\A)=0$, for  $q\le k$, as stated, via \eqref{eq:twistedf}.
\end{proof}

\subsection{}
\label{ss34}
Our theorem above complements a similar result obtained by Libgober, who proved in
\cite{Li2}, with a different method, that the non-divisibility conditions for all 
$X\in \L(\A)$, dense, with rank between $2$ and $k+1$, and contained in some $K\in \A$,
imply the same conclusion. Either vanishing criterion may be used to deduce the
following consequence, that led us to the formula from 
Theorem \ref{thm:main} \eqref{a4}.

\begin{corollary}
\label{cor:conj}
Let $\A$ be an arbitrary subarrangement, with $n$ hyperplanes and of rank $\ge 3$, 
of a Coxeter
arrangement. If $d \mid n$ and $d \not\in \{1,2,3,4,5\}$,
then $b_{1d}(\A)=0$.
\end{corollary}

\begin{proof}
We know that $\A \subset T$, where $T$ is a full Coxeter arrangement 
and $\rk (T)\ge 3$. Pick any
rank two element $X\in \L(\A)$. Plainly, $\A_X \subset T_X$.
Inspecting the tables from \cite{OT}, we conclude that $m_X\le 5$. Therefore,
the ${\bf m}_2$--list of $\A$ defined in \eqref{eq:mlist} is contained in
$\{ 3,4,5 \}$. Our assertion becomes then a direct consequence of 
Theorem \ref{thm:negvanish}.
\end{proof}
 
\section{Mod $p$ Aomoto complexes of graphic arrangements ($p\ne 3$)}
\label{sect4}

\subsection{}
\label{ss41}
We will use the following terminology and notation.
Denote by $[\ell]$ the set of points  $\{1,\dots,l\}$. We say that
$\Gamma$ is a graph in $[\ell]$ if the
set of edges of $\Gamma$ decomposes, $\E(\Gamma)=\E_{1}(\Gamma)\sqcup \E_{2}(\Gamma)$, 
where 
$\E_{1}(\Gamma)\subset [\ell]$ is the set of loops  and $\E_{2}(\Gamma)$, 
the set of signed edges, consists of elements 
of the form $ij^{\epsilon}$, with $\{i\ne j\} \subset [\ell]$ and $\epsilon \in \{\pm 1\}$.
 
\begin{definition}
\label{unsigned}
If $\Gamma$ is a graph in $[\ell]$, we denote by $\overline{\Gamma}$ the ordinary 
simplicial graph with set of edges 
$\E(\overline{\Gamma})=\{ij:= \{i\ne j\} \;|\; \exists\; \epsilon$
such that $ij^{\epsilon} \in \E_{2}(\Gamma)\}$.
We also denote by $\V(\Gamma)= \V(\overline{\G}):=
\{i\in [\ell]\;|\; \exists\; e \in \E(\overline{\Gamma})$ such that $i \in e\}$,
the set of vertices of $\Gamma$ ($\overline{\Gamma}$). 
\end{definition}

Here is the definition of the arrangement associated to 
a graph.

\begin{definition}
\label{def:agraph}
Let $\Gamma$ be a graph in $[\ell]$.
We denote by $\mathcal{A}(\Gamma)$ the arrangement in $\C^l$, with hyperplanes
given by the equations 
$x_{i}+ \epsilon x_{j}=0$, for each signed edge $ij^{\epsilon}\in \E_{2}(\Gamma)$, 
and $x_{i}=0$, 
for each loop $i \in \E_{1}(\Gamma)$.
\end{definition}

\begin{example}Complete graphs. 

\noindent \textit{(i)} If $\Gamma$ is the complete unsigned graph on $l$ vertices, 
then $\mathcal{A}(\Gamma)$ is the braid arrangement of rank $l-1$, with
defining equation
$\prod_{1\le i<j\le l} (x_i-x_j)=0$.

\noindent \textit{(ii)} If $\Gamma$ is the complete signed graph on $l$ vertices, 
then $\mathcal{A}(\Gamma)$ is the arrangement of hyperplanes corresponding 
to the Coxeter group $D_{l}$, with  
defining equation
$\prod_{1\le i<j\le l} (x_i\pm x_j)=0$.

\noindent \textit{(iii)} If in addition to that the graph has a loop at each vertex, 
then we get the arrangement corresponding to the Coxeter group $B_{l}$, defined by
$\prod_{i=1}^l x_i\cdot \prod_{1\le i<j\le l} (x_i\pm x_j)=0$.
\end{example}

\subsection{Rank 2 elements in a graphic arrangement}
\label{ss42}
 In what follows we will refer mainly to graphic arrangements, so it will be convenient 
 to use the label $\Gamma$ for objects associated to the arrangement $\mathcal{A}(\Gamma)$; 
 for instance, the lattice $\mathcal{L}(\mathcal{A}(\Gamma))$ is denoted 
 simply by $\mathcal{L}(\Gamma)$, and so on.
 
 \vskip 1in
 \begin{pspicture}(0,1.5)
 \label{m=2}
 \pscircle*(0,1){.07}
 \pscircle*(0,3){.07}
 \rput(0,0.5){\footnotesize $(1)$}
 \rput(0.2,2){\footnotesize $\pm$}
 \psline[linewidth=.5pt](0,1)(0,3)
 \rput(0.2,1){\footnotesize $i$}
 \rput(0.2,3){\footnotesize $j$}
 
 \pscircle*(2,1){.07}
 \pscircle*(2,3){.07}
 \rput(2.2,2){\footnotesize $\epsilon$}
 \psline[linewidth=.5pt](2,1)(2,3)
 \rput(2.2,1){\footnotesize $i$}
 \rput(2.2,3){\footnotesize $j$}
 \pscircle*(3,1){.07}
 \pscircle*(3,3){.07}
 \rput(2.5,0.5){\footnotesize $(2)$}
 \rput(3.2,2){\footnotesize $\epsilon'$}
 \psline[linewidth=.5pt](3,1)(3,3)
 \rput(3.2,1){\footnotesize $k$}
 \rput(3.2,3){\footnotesize $l$}
 
 \pscircle*(5,3){.07}
 \pscircle*(6,1){.07}
 \pscircle*(7,3){.07}
 \psline[linewidth=.5pt](5,3)(6,1)(7,3)
 \rput(6.2,1){\footnotesize $j$}
 \rput(4.8,3){\footnotesize $i$}
 \rput(7.2,3){\footnotesize $k$}
 \rput(6.7,2){\footnotesize $\epsilon'$}
 \rput(5.3,2){\footnotesize $\epsilon$}
 \rput(6,0.5){\footnotesize $(3)$}

 \pscircle*(9,1){.07}
 \pscircle*(9,3){.07}
 \pscircle(9,3){.4}
 \psline[linewidth=.5pt](9,1)(9,3)
 \rput(9,0.5){\footnotesize $(4)$}
 \rput(9.2,2){\footnotesize $\epsilon$}
 \rput(9.2,3){\footnotesize $i$}
 \rput(9.2,1){\footnotesize $j$}
 
 \pscircle*(11,1){.07}
 \pscircle*(11,3){.07}
 \pscircle*(12,2){.07}
 \pscircle(12,2){.4}
 \psline[linewidth=.5pt](11,1)(11,3)
 \rput(11.6,0.5){\footnotesize $(5)$}
 \rput(11.2,2){\footnotesize $\epsilon$}
 \rput(11.2,3){\footnotesize $i$}
 \rput(11.2,1){\footnotesize $j$}
 \rput(12.2,2){\footnotesize $k$}
 
 \pscircle*(14,1.2){.07}
 \pscircle(14,3){.4}
 \pscircle*(14,3){.07}
 \pscircle(14,1.2){.4}
 \rput(14.2,3){\footnotesize $i$}
 \rput(14.2,1.2){\footnotesize $j$}
 \rput(14,0.5){\footnotesize $(6)$}
 
 \end{pspicture}
 \vskip .1in
 
 \begin{figure}[ht]
 \caption{\textsf{Pairs of edges}}
 \label{fig:graph2}
 \end{figure}

 \vskip 0.5in
 \begin{pspicture}(0,1.5)
 \label{m=3,4}
 \pscircle*(0.4,1){.07}
 \pscircle*(2.4,1){.07}
 \pscircle(0.4,1){.4}
 \pscircle(2.4,1){.4}
 \rput(1.4,0.5){\footnotesize $(1)$}
 \rput(1.4,1.2){\footnotesize $\epsilon$}
 \psline[linewidth=.5pt](.4,1)(2.4,1)
 \rput(0.4,1.3){\footnotesize $i$}
 \rput(2.4,1.3){\footnotesize $j$}
 
 \pscircle*(4.4,1){.07}
 \pscircle*(6.4,1){.07}
 \pscircle(4.4,1){.4}
 \pscircle(4.4,1){.4}
 \rput(5.4,0.5){\footnotesize $(2)$}
 \rput(5.4,1.2){\footnotesize $\pm$}
 \psline[linewidth=.5pt](4.4,1)(6.4,1)
 \rput(4.4,1.3){\footnotesize $i$}
 \rput(6.4,1.3){\footnotesize $j$}

 \pscircle*(7.9,1){.07}
 \pscircle*(9.9,1){.07}
 \pscircle*(8.9,2.4){.07}
 \psline[linewidth=.5pt](8.9,2.4)(7.9,1)(9.9,1)(8.9,2.4)
 \rput(8.9,0.5){\footnotesize $(3)$}
 \rput(8.9,1.2){\footnotesize $\epsilon$}
 \rput(8.3,1.8){\footnotesize $\epsilon'$}
 \rput(9.7,1.8){\footnotesize $\epsilon''$}
 \rput(7.7,1.2){\footnotesize $i$}
 \rput(10.1,1.2){\footnotesize $j$}
 \rput(9.1,2.6){\footnotesize $k$}
 
 \pscircle*(11.9,1){.07}
 \pscircle*(13.9,1){.07}
 \pscircle(11.9,1){.4}
 \pscircle(13.9,1){.4}
 \rput(12.9,1.2){\footnotesize $\pm$}
 \psline[linewidth=.5pt](11.9,1)(13.9,1)
 \rput(11.7,1.2){\footnotesize $i$}
 \rput(14,1.3){\footnotesize $j$}
 \rput(12.9,0.5){\footnotesize $(4)$}
 
 \end{pspicture}
 \vskip .1in
 
 \begin{figure}[ht]
 \caption{\textsf{Dense elements}}
 \label{fig:graph34}
 \end{figure}

 For reasons that will become clear from subsection \S\ref{ss44} on, we 
 draw up a complete inventory of 
 rank $2$ elements $X \in \mathcal{L}(\Gamma)$, by representing the subgraphs corresponding 
 to the associated subarrangements, $\mathcal{A}_{X}(\Gamma)$. See figures
 \ref{fig:graph2} and \ref{fig:graph34}.
 
 \begin{remark}
 Recall from \S\ref{ss31} that $m_{X}$ denotes the number of hyperplanes in 
 the subarrangement $\mathcal{A}_{X}$,
 for $X\in \mathcal {L}(\mathcal{A})$. In Figure \ref{fig:graph2}, $m_X=2$, while
 $m_X =3$ or $4$, in Figure \ref{fig:graph34}.
 In Figure \ref{fig:graph2}, the configuration $(3)$ means that
 $ik^{- \epsilon \epsilon'} \notin \E_{2}(\Gamma)$.
 In Figure \ref{fig:graph2}$(4)$, $ij$ is a
 {\em simple} edge of $\overline{\G}$ (identified with the corresponding edge, 
 $ij^{\epsilon}$, of $\G$), that is, $ij^{-\epsilon}\notin \E_2(\G)$.
 In Figure \ref{fig:graph34}$(2)$, $ij$ is a
 {\em double} edge of $\overline{\G}$ (identified with the corresponding pair of edges in $\G$, 
 $ij^{\pm}$).
 In Figure \ref{fig:graph34}$(3)$, the signs on the edges must be such that 
 $\epsilon \epsilon' \epsilon''=-1$. Such
 a triangle is called \textit{negative} (otherwise the triangle is called \textit{positive}).
 \end{remark}

\subsection{Weighted graphs}
\label{ss43}

 An element $\eta \in A^{1}_{\k}(\Gamma)$, $\k$  a field, may be viewed as
 a collection of \textit{weights}, that is, a set of coefficients, 
 $\eta_{k} \in\k$, one
 for each $k \in \E_{1}(\Gamma)$, and $\eta _{ij}^{\epsilon}\in \k$, one for 
 each $ij^{\epsilon} \in \E_{2}(\Gamma)$.
 If $\k= \mathbb{F}_{p}$, we will abbreviate $ \mathbb{F}_{p}$ by $p$,
 when referring to the coefficient field; for instance, $A_{p}^{1}(\Gamma) :=
 A_{\mathbb{F}_{p}}^{1}(\Gamma)$.
 
\begin{remark}
\label{constant}
Denote by $Z_{p}(\Gamma)$ the set of $1$--cocycles in $(A_{p}^{*}(\Gamma),\mu_{1})$
(see \eqref{eq:aok}). Then $\beta_{p}(\Gamma)=0$ if and only if 
the weights of $\eta$ are constant on $\E(\G)$,
for any $\eta \in Z_{p}(\Gamma)$.
\end{remark}
 
The following well-known result will be extensively used 
in computing $\beta_{p}(\Gamma)$, for $p$ a prime.
  
\begin{lemma}
\label{rules}
Let $\mathcal{A}$ be an arbitrary central arrangement, $p$ be a prime. If $\eta \in A_{p}^{1}
(\mathcal{A})$, $\eta =\sum_{H\in \A}\eta_H a_H$,
then $\eta \omega_{1}=0$ if and only if one has 
\begin{equation}
\label{div}
\Sigma_X \eta:=\sum_{H \supset X} \eta_{H}=0, \;\;if \;\;p \mid m_{X},
\end{equation}
or 
\begin{equation}
\label{ndiv}
\eta_{H}= \eta_{K},\;\;\forall\; H\ne K \;\;\in \mathcal{A}_{X}, \;\;if\; p\;\; \nmid m_{X}\,,
\end{equation}
for every rank $2$ element 
$X \in \mathcal{L}_{2}(\mathcal{A})$.
\end{lemma}

\begin{proof}
See for instance \cite[Lemma 3.3]{LY}.
\end{proof}

\subsection{Graphic arrangements at primes different from $3$}
\label{ss44}

We will need to compute the numbers $\beta_{p}(\Gamma)$, for arbitrary
$\G$ and $p$, when $\rk\; \mathcal{A}(\Gamma) >2$. We end this section by
showing that these numbers are zero, for $p\ne 3$.

\begin{lemma}
\label{pnot234}
If $p \neq 2,\;3$, then $\beta_{p}(\Gamma)=0$.
\end{lemma}

\begin{proof}
Let $H \neq K$ be arbitrary hyperplanes in $\mathcal{A}(\Gamma)$. Set 
$X=X(H,K):=H \cap K \in \mathcal{L}_{2}(\Gamma)$. Consider $\eta \in Z_{p}(\Gamma)$, 
$\eta=\sum_{H \in \mathcal{A}(\Gamma)} \eta_{H}a_{H}$.
By inspecting Figures \ref{fig:graph2} and \ref{fig:graph34} 
from subsection \ref{ss42}, we see
that the condition $p\nmid m_{X}$ from \eqref{ndiv} is satisfied, 
so $\eta_{H}= \eta_{K}$, as
needed (see Remark \ref{constant}).
\end{proof}

The same argument actually proves the following analog of Theorem \ref{thm:negvanish}.

\begin{prop}
\label{prop=betav}
Let $\A$ be an arbitrary central arrangement. 
If a prime $p$ does not divide $m$, for any $m\in {\bf m}_{2}(\A)$, then
$\beta_{1p}(\A)=0$.
\end{prop}

\begin{corollary}
\label{cor=betaconj}
Let $\A$ be an arbitrary subarrangement, of rank $\ge 3$, 
of an arbitrary Coxeter arrangement. Then $\beta_{1p}(\A)=0$, for
$p\notin \{ 2, 3, 5\}$.
\end{corollary}

\begin{proof}
Recall from the proof of Corollary \ref{cor:conj} that
${\bf m}_{2}(\A) \subset \{ 3, 4, 5\}$.
\end{proof}

\begin{prop}
\label{p2}
Assume $\rk\; \mathcal{A}(\Gamma) >2$. Then $\beta_{2}(\Gamma)=0$.
\end{prop}

\begin{proof}
Consider an arbitrary element $\eta \in Z_{2}(\Gamma)$.
We have to show that $\eta_{H}= \eta_{K}$, $\forall H\ne K\in \A(\G)$. 
Set $X=H\cap K\in \L_2(\G)$.
If $m_{X} \in \{2,3\}$, then we are done, by resorting
to Lemma \ref{rules}.

Otherwise, $m_{X}=4$, that is,
the subarrangement $\mathcal{A}_{X}(\G)$ is given by a subgraph of the type
depicted in Figure \ref{fig:graph34}$(4)$, where say $i=1$ and $j=2$.

Then the weights of $\eta$ on  $\mathcal{A}_{X}(\G)$ must satisfy 
\begin{equation}
\label{suma}
\eta_{12}^{+}+\eta_{12}^{-}+\eta_{1}+\eta_{2}=0\, ,
\end{equation}
by Lemma \ref{rules}.
Since $\rk \mathcal{A}(\Gamma) >2$, there must be an edge $e$ 
(of weight say $a$) in $\E(\Gamma)$,
corresponding to a hyperplane that does not contain $X$. 

Two cases may occur:

{\em Case} $(a)$ There is an edge $e \in \E_2(\Gamma)$, different from $12^{\pm}$.

{\em Subcase} $(a.1)$ Both endpoints of $e$ are different from $1$ and $2$. In this case,
figures \ref{fig:graph2}(2) and \ref{fig:graph2}(5) imply, 
via Lemma \ref{rules}, that $\eta$
has constant weight, equal to $a$, on $\mathcal{A}_{X}(\G)$. In particular,
$\eta_{H}=\eta_{K}$, as asserted. 

{\em Subcase} $(a.2)$ Otherwise, we may assume $e=13^{\epsilon}\in \E_2(\G)$. 
Then $\eta_{2}=a$ 
(see figure \ref{fig:graph2}(5) and Lemma \ref{rules}). Moreover,
$\eta_{12}^{-}=\eta_{12}^{+}=a$, as follows from 
figure \ref{fig:graph2}(3) or figure \ref{fig:graph34}(3), again by Lemma \ref{rules}. 
We infer then from \eqref{suma} that $\eta$
has constant weight on $\mathcal{A}_{X}(\G)$, and we are done.

{\em Case} $(b)$ There are no other edges in $\E_2(\Gamma)$, except $12^{\pm}$, but
there is a loop $e$ in $\E_1(\Gamma)$, at $k\ne 1, 2$. Then
$\eta_{12}^{\pm}=a$, and $\eta_1=\eta_2=a$, by Lemma \ref{rules}
(see figure \ref{fig:graph2}(5) and figure \ref{fig:graph2}(6) respectively).
\end{proof}

\section{Mod $3$ graphic Aomoto complexes}
\label{sect43}

We analyze now what happens at the prime $3$.

\begin{prop}
\label{p3}
Assume  $\rk\; \mathcal{A}(\Gamma) >2$. 
\begin{enumerate}
\item \label{3notexc}
If $\beta_{3}(\Gamma)\neq 0$, 
then $\Gamma$ must be one of the 
graphs from Figures \ref{5graphs} and \ref{2moregraphs}.

\item \label{3exc}
If $\G$ is exceptional, then $\beta_3(\G)=1$.
\end{enumerate}
\end{prop}

\subsection{Preliminary lemmas}
\label{3lem}
The proof of Proposition \ref{p3} will occupy the rest of this section,
where the coefficient field is understood to be $\mathbb{F}_{3}$.

\vskip 1.5in
\begin{pspicture}(0,1.5)
\label{exceptions}

\pscircle*(1,1){.07}
\pscircle*(3,1){.07}
\pscircle*(2,2.4){.07}
\psline[linewidth=.5pt](2,2.4)(1,1)(3,1)(2,2.4)
\rput(2,0.5){\footnotesize $(1)$}
\rput(2,1.2){\footnotesize $\pm$}
\rput(1.4,1.8){\footnotesize $\pm$}
\rput(2.7,1.8){\footnotesize $\pm$}

\pscircle*(5,1){.07}
\pscircle*(7,1){.07}
\pscircle*(6,2.4){.07}
\pscircle(5,1){.4}
\pscircle(7,1){.4}
\pscircle(6,2.4){.4}
\psline[linewidth=.5pt](6,2.4)(5,1)(7,1)(6,2.4)
\rput(6,0.5){\footnotesize $(2)$}
\rput(6,1.2){\footnotesize $-\epsilon \epsilon'$}
\rput(5.2,1.8){\footnotesize $\epsilon'$}
\rput(6.8,1.8){\footnotesize $\epsilon$}

\pscircle*(9,1){.07}
\pscircle*(11,1){.07}
\pscircle(10,2.4){.4}
\pscircle*(10,2.4){.07}
\psline[linewidth=.5pt](10,2.4)(9,1)(11,1)(10,2.4)
\rput(10,0.5){\footnotesize $(3)$}
\rput(10,1.2){\footnotesize $\epsilon$}
\rput(9.4,1.8){\footnotesize $\pm$}
\rput(10.7,1.8){\footnotesize $\pm$}

\end{pspicture}

\begin{figure}[ht]
\caption{\textsf{Exceptional graphs}}
\label{5graphs}
\end{figure}

\vskip .5in
\begin{pspicture}(0,1.5)
\label{more exceptions}

\pscircle*(2,-3){.07}
\pscircle*(5,-3){.07}
\pscircle*(3.5,-2){.07}
\pscircle*(3.5,0){.07}
\psline[linewidth=.5pt](2,-3)(3.5,0)(5,-3)(2,-3)(3.5,-2)(5,-3)
\psline[linewidth=.5pt](3.5,0)(3.5,-2)
\rput(3.5,-3.8){\footnotesize $(4)$}
\rput(3.5,-3.2){\footnotesize $-\epsilon' \epsilon''$}
\rput(2.4,-1.5){\footnotesize $-\epsilon \epsilon'$}
\rput(4.7,-1.5){\footnotesize $-\epsilon \epsilon''$}
\rput(3.7,-1){\footnotesize $\epsilon$}
\rput(4.3,-2.2){\footnotesize $\epsilon'' $}
\rput(2.8,-2.2){\footnotesize $\epsilon'$}

\pscircle*(8,-3){.07}
\pscircle*(11,-3){.07}
\pscircle*(9.5,-2){.07}
\pscircle*(9.5,0){.07}
\psline[linewidth=.5pt](8,-3)(9.5,0)(11,-3)(8,-3)(9.5,-2)(11,-3)
\psline[linewidth=.5pt](9.5,0)(9.5,-2)
\rput(9.5,-3.8){\footnotesize $(5)$}
\rput(9.5,-3.2){\footnotesize $\pm$}
\rput(8.4,-1.5){\footnotesize $\pm$}
\rput(10.6,-1.5){\footnotesize $\pm$}
\rput(9.7,-1){\footnotesize $\pm$}
\rput(10.2,-2.2){\footnotesize $\pm $}
\rput(8.8,-2.2){\footnotesize $\pm$}

\end{pspicture}
\vskip 1.5in

\begin{figure}[ht]
\caption{\textsf{More exceptional graphs}}
\label{2moregraphs}
\end{figure}

\begin{lemma}
\label{triunghi}
Let $\Gamma'\subset \Gamma$ be a subgraph such that 
$\overline{\G'}$ is a triangle. Assume that $\E_2(\G')$ contains
a simple edge of $\overline{\G}$, and a
double edge of $\overline{\G}$. 
Assume also that $\G$ has no loops at the vertices of the triangle.
If $\eta \in Z(\G)$,
then the weights of $\eta$ are constant,
on all edges of $\G'$.
\end{lemma}

\begin{proof}
Let the subgraph be as in the picture below. Here the edge $13$ is double
($13^{\pm} \in \E_2(\G')$),
the edge $12$ is simple ($12^{\epsilon} \in \E_2(\G')$, $12^{-\epsilon} \notin \E_2(\G)$), 
and $23^{\epsilon'}$  is one of the
(at most two) edges from $\E_2(\G')$ corresponding to $23 \in \E(\overline{\G'})$.
Denote $\eta_{13}^{+}$ by $a$. We have to show that
$\eta_{13}^{-}=\eta_{12}^{\epsilon}= \eta_{23}^{\epsilon'}=a$.

\vskip .8in
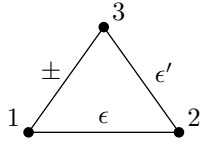
\begin{pspicture}(5,1.5)
\pscircle*(6,1){.07}
\pscircle*(8,1){.07}
\pscircle*(7,2.4){.07}
\psline[linewidth=.5pt](7,2.4)(6,1)(8,1)(7,2.4)
\rput(7,1.2){\footnotesize $\epsilon$}
\rput(6.3,1.8){\footnotesize $\pm$}
\rput(7.8,1.8){\footnotesize $\epsilon'$}
\rput(5.8,1.2){\footnotesize $1$}
\rput(8.2,1.2){\footnotesize $2$}
\rput(7.2,2.6){\footnotesize $3$}

\end{pspicture}
\vskip .01in
 
 Since there are no $\G$--loops in $[3]$, we infer from 
 figure \ref{fig:graph2}(1)
 and Lemma \ref{rules} that
 $\eta_{13}^{+}=\eta_{13}^{-}=a$. 

 Set $\epsilon''=\epsilon \epsilon'$. Then  
 $\eta_{23}^{\epsilon'}= \eta_{13}^{\epsilon''}= a$, since
 ${12}^{-\epsilon}\notin \E_2(\G)$
 (see figure \ref{fig:graph2}(3) and \eqref{ndiv}). Next, 
 we obtain from figure \ref{fig:graph34}(3) and 
 \eqref{div} that
 $\eta_{12}^{\epsilon}+ \eta_{13}^{-\epsilon''} + \eta_{23}^{\epsilon'}= 0$. Therefore,
 $\eta_{12}^{\epsilon}+ 2a=0$, whence $\eta_{12}^{\epsilon}=a$.
 Consequently, all weights of $\eta$ from  the triangle above are equal to $a$.
\end{proof}

The following definition will be convenient for our purposes:
the {\em full subgraph} $\G'$ of $\G$, determined by $\V'\subset \V(\G)$, has edges 
$\E(\G')=\E_2(\G'):= \{ ij^{\epsilon}\in \E_2(\G) \, \mid \, i,j\in \V' \}$.

\begin{lemma}
\label{tetraedru}
 Let $\Gamma$ be a graph whose associated unsigned graph, $\overline{\G}$, is
 complete on $4$ vertices. If $\eta\in Z(\G)$ has 
 constant weight on $\E_{2}(\Gamma')$,
 where $\G'$ is a full subgraph of $\G$ on $3$ vertices, then $\eta$
 has constant weight on $\E_{2}(\Gamma)$.
\end{lemma}

\begin{proof}
 Set $\V(\Gamma')=[3]\subset [4]=\V(\Gamma)$. 
 We know that $\eta$ has weight
 $a$, on $\E_{2}(\Gamma')$. Pick any edge 
 $e=ij^{\epsilon}\in \E_{2}(\Gamma)\setminus \E_{2}(\Gamma')$.
 Clearly, $\abs{\{i, j\}\cap [3]}=1$, since $\G'$ is the full subgraph of $\G$ 
 determined by $[3]$.
 Hence, we may find
 another edge, $f=kl^{\epsilon'}\in \E_{2}(\Gamma')$, such that 
 $\{ i, j\}\cap \{k, l\}=\emptyset$. Figure \ref{fig:graph2}(2) and \eqref{ndiv}
 together imply that $\eta_{ij}^{\epsilon}=\eta_{kl}^{\epsilon'}=a$.
\end{proof}

\begin{lemma}
\label{loop}
 Let $\Gamma$ be a graph whose associated unsigned graph, $\overline{\G}$, is
 complete on $4$ vertices. If $\E_1(\G)\ne \emptyset$, then 
 the weights of $\eta$ on $\Gamma$ are constant, for any $\eta\in Z(\G)$.
\end{lemma}

\begin{proof}
 Let $i\in \E_1(\G)$ be an arbitrary loop, with weight $a$. We have to show that
 $\eta$ has constant weight $a$ on $\E_2(\G)$. We may 
 assume that $\V(\Gamma)=[4]$, and $i=4$.
 (Indeed, if $i\notin \V(\G)$, then figure \ref{fig:graph2}(5) and Lemma \ref{rules}
 give the desired conclusion.)
 Then $\eta_{ij}^{\epsilon}=a$, for any edge $ij^{\epsilon}$ of the full subgraph
 of $\G$ determined by $[3]$
 (use figure \ref{fig:graph2}(5) and \eqref{ndiv}).
 Lemma \ref{tetraedru} yields then the desired conclusion.
\end{proof}

\subsection{}
\label{prel}

We begin the proof of Proposition \ref{p3}\eqref{3notexc} by a few preliminary
remarks.

\begin{remark}
\label{rk:3cases}
The assumption $\beta_3(\G)\ne0$ guarantees the existence of $\eta\in Z_3(\G)$ with
the property that the weights of $\eta$ are not constant on $\A_X(\G)$, for some
$X\in \L_2(\G)$. By Lemma \ref{rules}\eqref{ndiv}, this forces 
$m_{X}=3$. In other words, the subarrangement $\mathcal{A}_X(\G)$ is represented by
one of the first three graphs from Figure \ref{fig:graph34}. 
So, there are three cases to be examined.
\end{remark}

 \begin{remark}
 \label{+}
 For each of the above configurations, the fact that two
 out of the three weights of $\eta$ on $\A_X(\G)$ are equal 
 is equivalent  to the fact that  
 $\eta$ has constant weight on $\A_X(\G)$
 (use \eqref{div} and remember that we are working modulo $3$).
 \end{remark}

 \begin{remark}
 \label{rk}
 Due to our assumption on $\rk \mathcal{A}(\Gamma)$, there must be 
 an edge $e\in \E(\Gamma)$, different from those of $\mathcal{A}_{X}(\G)$.
 \end{remark}
 
 \subsection{Proof of Proposition \ref{p3}\eqref{3notexc}}
 \label{proofnotexc}
 
 We proceed to the analysis of the 3 abovementioned cases. Whenever
 possible without creating any ambiguity,
 we will omit the non-relevant signs of  edges from $\E_2(\G)$, to avoid making
 the exposition too heavy.

 {\em Case} \textbf{(a)}: Suppose $\A_{X}(\G)$ corresponds to a subgraph in $\Gamma$ 
 of the type described in Figure
 \ref{fig:graph34}(3), with vertices labeled $i=1$,$j=2$,$k=3$.
 We know that $\eta_{12}+\eta_{23}+\eta_{13}=0$, from Lemma \ref{rules}\eqref{div}.
 
 $(a.0)$  We may assume in case  \textbf{(a)}
 that there is no edge in $\G$ of the form $e=ij$, with
 $\{ i, j\}\cap [3]=\emptyset$. Indeed, otherwise  figure \ref{fig:graph2}(2)
 and \eqref{ndiv} would imply that all weights of $\eta$ on $\A_X(\G)$ are equal 
 to the weight of $e$, in contradiction with Remark \ref{rk:3cases}.

 Our discussion splits now, according to the number of vertices of $\G$: either
 $\abs{\V(\G)}>3$, or  $\abs{\V(\G)}=3$.
 
 {\em Case} $(a.1)$:  $\abs{\V(\G)}>3$. 
 We first claim that $ij\in \E(\overline{\G})$, for every 
 vertex $j$ of $\G$, $j\notin [3]$, and for all $i\in [3]$. 
 
 Indeed, denoting $j$ by $4$,
 we may resort to $(a.0)$ 
 to assume that say $14$ is 
 an edge of $\G$, with weight $a$. Then
 $\eta_{23}=a$ (by Lemma \ref{rules}, applied to figure \ref{fig:graph2}(2)).
 
 If there is no edge in $\G$ connecting the vertices 2 and 4, or 3 and 4, 
 we may apply Lemma \ref{rules} to figure \ref{fig:graph2}(3) to deduce that 
 $a=\eta_{12}$
 (respectively $a=\eta_{13}$). Hence, $\eta$ must be constant on $\A_{X}(\G)$
 (see Remark \ref{+}), which contradicts Remark \ref{rk:3cases}. The claim
 is thus verified.

 Again, there are two possibilities: either  $\abs{\V(\G)}>4$, or  $\abs{\V(\G)}=4$.

 {\em Subcase} $(a.1.1)$: There are another vertices, say $4$ and $5$, of $\G$. 
 Due to the previous claim, $ij\in \E(\overline{\G})$, for all $i\in [3]$
 and $j=4,5$. It follows that $\eta_{12}=\eta_{35}$,  $\eta_{13}=\eta_{24}$,
 and $\eta_{24}=\eta_{35}$, see figure \ref{fig:graph2}(2) and \eqref{ndiv}.
 Therefore, $\eta_{12}=\eta_{13}$, contradicting again Remark \ref{rk:3cases}, via
 Remark \ref{+}.

 {\em Subcase} $(a.1.2)$: $\V(\G)=[4]$. We already know that $\overline{\G}$ is a
 complete graph.
  
 If there exists a loop in $\Gamma$, we obtain a contradiction 
 by applying Lemma \ref{loop}.
 So, there are no loops in $\Gamma$. Now, if $\Gamma$ contains a full subgraph 
 on $3$ vertices, having both
 simple and  double edges, we may invoke lemmas \ref{triunghi} and \ref{tetraedru} to
 infer that $\eta$ has constant weight on $\G$, which leads to the same contradiction 
 as before. If not, it follows that $\Gamma$ must be one of the graphs from
 Figure \ref{2moregraphs}.
 
 Indeed, this is clear if all edges of $\overline{\G}$ are double. 
 Otherwise, they must be all
 simple. Now, if there is a positive triangle in $\G$, then $\eta$ must have constant
 weight on it (see figure \ref{fig:graph2}(3)). Again, Lemma \ref{tetraedru} leads to
 a contradiction.
 
 This completes the discussion of Case $(a.1)$.

 {\em Case} $(a.2)$: $\V(\Gamma)=[3]$. In this case, we may suppose $\E_1(\G)\subset [3]$
 (otherwise, the equations provided by figure \ref{fig:graph2}(5) would force $\eta$ to
 have constant weight on $\E_2(\G)$, in particular on $\A_X(\G)$). 
 In what follows, the discussion naturally splits according
 to the number of loops in $\G$.

 {\em Subcase} $(a.2.0)$: There are no loops in $\Gamma$. 
 By virtue of Lemma \ref{triunghi}, all edges must be
 double (see Remark \ref{rk}). Thus, $\Gamma=D_{3}$, the first graph from Figure
 \ref{5graphs}.

 {\em Subcase} $(a.2.1)$:  $\abs{\E_{1}(\Gamma)}=1$. 
 Let $1$ be the unique loop, with weight $a$. Then
 $a=\eta_{23}$ (by Lemma \ref{rules} and figure \ref{fig:graph2}(5)). 
 At this point, two possibilities may occur.
 
 {\em Subcase} $(a.2.1')$:
 One of the  edges $12$ or $13$ is simple. In this situation, 
 we may apply Lemma \ref{rules} to
 figure \ref{fig:graph2}(4), deducing that either $\eta_{12}=a$ or $\eta_{13}=a$,
 which contradicts Remark \ref{rk:3cases} (see Remark \ref{+}).
 
 {\em Subcase} $(a.2.1'')$: Otherwise, both edges $12$ and $13$ are double. 
 When all edges are double,  $\G$ is entirely determined; a routine application
 of Lemma \ref{rules} shows then that
 $\beta_3(\G)=0$. When the edge $23$ is simple, 
 we obtain the  graph from Figure \ref{5graphs}(3).

 {\em Subcase} $(a.2.2)$:  $\abs{\E_{1}(\Gamma)}=2$, i.e., $\E_{1}(\Gamma)$ is say $\{1, 2\}$.
 
 {\em Subcase} $(a.2.2')$: The edge $12$ is double. Then it follows from Lemma \ref{rules}
 \eqref{ndiv} that all $4$ edges of the configuration from
 Figure \ref{fig:graph34}(4) (where $ij=12$)
 have the same weight, say $a$. 
 
 If one of the other edges, say $23$, is simple,
 Lemma \ref{rules}\eqref{ndiv} may be applied to figure
 \ref{fig:graph2}(4), to infer that $\eta_{23}=a$. By Remark \ref{+}, this 
 contradicts Remark \ref{rk:3cases}.
 
 Finally, if all edges are double, a straightforward computation shows that
 $\beta_3(\G)=0$, like in  subcase $(a.2.1'')$.

 {\em Subcase} $(a.2.2'')$: The edge $12$ is simple. This implies that 
 $\eta_1+\eta_2 +\eta_{12}=0$ (see figure \ref{fig:graph34}(1)).
 If the edge $13$ is also simple, we obtain $\eta_1=\eta_{13}$ 
 (see figure \ref{fig:graph2}(4)).
 We also get, by using figure \ref{fig:graph2}(5), that 
 $\eta_2=\eta_{13}$. Putting these
 facts together, we deduce that $\eta_{12}=\eta_{13}$, a contradiction.
 If the edge $13$ is double, then 
 $\eta_1+\eta_{13}^+ +\eta_{13}^- =0$ (see figure \ref{fig:graph34}(2)),
 and $\eta_{13}^{\pm} =\eta_{2}$ (see figure \ref{fig:graph2}(5)). Hence, the weights
 $\eta_1$, $\eta_2$, $\eta_{12}$ and  $\eta_{13}^{\pm}$ are all equal. In particular, 
 $\eta_{12}= \eta_{13}^{\epsilon'}$, a contradiction.
 
 {\em Subcase} $(a.2.3)$:  $\E_{1}(\Gamma)=[3]$.
 
 {\em Subcase} $(a.2.3')$: There is a simple edge, say $12$, 
 and a double edge, say $13$.
 In this case, we have: $\eta_1=\eta_3=\eta_{13}^{\pm}$ 
 (see figure \ref{fig:graph34}(4)),
 and $\eta_3=\eta_{12}$ (see figure \ref{fig:graph2}(5)). These facts yield
 $\eta_{12}=\eta_{13}^{\epsilon'}$, a contradiction, as before.
 
 {\em Subcase} $(a.2.3'')$: Either all edges are simple, i.e., $\G$ is the graph from
 Figure \ref{5graphs}(2), or all edges are double, and then  
 it is easy to see that $\beta_3(\G)=0$.
 
 The analysis of case $\textbf{(a)}$ is thus complete.

 In the remaining two cases, $\A_X(\G)=\A(\G')$, where $\G'$ is a subgraph with
 shape described in figure \ref{fig:graph34}(1)--(2), with say $ij=12$. We begin by
 two remarks, valid in both these cases.
 
 $(bc.1)$ We may assume that there is no edge $ij$ in $\G$ disjoint from $12$.
 Indeed, otherwise figures \ref{fig:graph2}(2) and \ref{fig:graph2}(5) 
 would imply, via Lemma \ref{rules},
 that all weights of $\eta$ on $\A_X(\G)$ are equal to the weight of $ij$, a
 contradiction.
 
 $(bc.2)$ We may also assume that $\E_2(\G)\ne \E_2(\G')$. If not, Remark
 \ref{rk} guarantees the existence of a loop of $\G$ away from $[2]$, 
 say $3$. Using this time
 figures \ref{fig:graph2}(5) and \ref{fig:graph2}(6), 
 we arrive again at a contradiction, as before.

 {\em Case} \textbf{(b)}: $\A_{X}(\G)$ corresponds to a subgraph in $\Gamma$ 
 of the type from figure \ref{fig:graph34}(1). We know from Lemma \ref{rules} that 
 $\eta_1+\eta_2+\eta_{12}=0$.
 
 It follows from $(bc.1)-(bc.2)$ above
 that we may suppose $13\in \E(\overline{\G})$. If 
 $23\notin \E(\overline{\G})$, we infer 
 from lemma \ref{rules} that
 $\eta_{12}=\eta_{13}$ and $\eta_{2}=\eta_{13}$
 (see figure \ref{fig:graph2}, (3) and (5) respectively), thus contradicting 
 Remark \ref{rk:3cases}, via Remark \ref{+}. It follows that 
 $13^{\epsilon'}, 23^{\epsilon''}\in \E_2(\G)$, for some signs,
 $\epsilon'$ and $\epsilon''$.
 
 {\em Subcase} $(b+)$: The triangle $123$ is positive. Then 
 $\eta_{13}^{\epsilon'}=\eta_{23}^{\epsilon''}$ (see figure \ref{fig:graph2}(3)).
 Moreover, $\eta_{1}=\eta_{23}^{\epsilon''}$ and $\eta_{2}=\eta_{13}^{\epsilon'}$
 (see figure \ref{fig:graph2}(5)). Hence, $\eta_{1}=\eta_{2}$, a contradiction again.
 
 {\em Subcase} $(b-)$: The triangle $123$ is negative. If the weights of $\eta$ on this
 triangle are not constant, we are back in case \textbf{(a)}, and  we are done. 
 Otherwise, denoting by $a$ their common value, we may use figure \ref{fig:graph2}(5) to
 deduce that $\eta$ must have constant weight $a$ on $\A_{X}(\G)$, which  
 contradicts our initial assumption from Remark \ref{rk:3cases}.
 
 The analysis of Case \textbf{(b)} is thus completed.
 
 {\em Case} \textbf{(c)}: $\A_{X}(\G)$ corresponds to a subgraph in $\Gamma$ 
 of the type from figure \ref{fig:graph34}(2). Lemma \ref{rules} implies that
 $\eta_1+ \eta_{12}^+ +\eta_{12}^-=0$.
 
 As before, we know that
 either $13$ or $23$ is an edge of $\G$, of weight say $a$. 
 If they do not both belong to $\E (\overline{\G})$,
 then figure \ref{fig:graph2}(3) forces $\eta_{12}^{\pm}=a$, a contradiction.
 Consequently, we may find a negative triangle in $\G$, with edges
 $13^{\epsilon'}$, $23^{\epsilon''}$ and $12^{\epsilon}$.
 Moreover, $\eta_1=\eta_{23}^{\epsilon''}$
 (see figure \ref{fig:graph2}(5)). 
 
 If $\eta$ has constant weight $a$ on this triangle, then $\eta_1=\eta_{12}^{\epsilon}=a$.
 Therefore, $\eta$ must also have
 constant weight $a$ on $\A_{X}(\G)$, by Remark \ref{+}, which is impossible.
 Otherwise, we are again back in case \textbf{(a)}, and  we are done. 
 
 This finishes the proof of Proposition \ref{p3}\eqref{3notexc}.
 
 \subsection{}
 \label{end5}
 
 Proposition \ref{p3}\eqref{3exc} will follow from the next two lemmas.

\begin{lemma}
\label{D3D4}
$\beta_{3}(D_{3})=\beta_{3}(D_{4})=1$.
\end{lemma}

\begin{proof}
Direct computation, using Lemma \ref{rules}.
\end{proof}

\begin{lemma}
\label{liso}
The exceptional graphic arrangements from Figures \ref{5graphs}(2)-(3) 
and \ref{2moregraphs}(4)
are lattice--isotopic to $D_3$.
\end{lemma}

\begin{proof}
We begin with the simplest case: the graph $\G$ from figure  \ref{2moregraphs}(4).
By a convenient change of signs of the variables from $\C^4$, we can transform
$\A(\G)$ into $A_3=D_3$. Similarly, we may assume that $\epsilon=\epsilon'=-1$, for 
the graphic arrangement $\A(\G)$ from figure \ref{5graphs}(2); by an obvious linear
change of coordinates, we can finally make  $\A(\G)$ projectively equivalent, hence
lattice--isotopic, to $D_3$.

By a preliminary change of signs, the last arrangement $\A(\G)$ 
(see figure \ref{5graphs}(3))
becomes defined by the equation $x_1(x_1\pm x_2)(x_1\pm x_3)(x_2-x_3)=0$.
Next, we make the change of variables $x_1=z_2+z_3$; $x_1+x_2=z_1+z_3$; $x_1+x_3=z_1+z_2$.
We arrive at a defining equation that corresponds to the value $t=-1$ in the 
family below (where $t\ne 1$)
\[
(z_1+z_2)(z_1+z_3)(z_2\pm z_3)[(z_1- z_2)+ t(z_2+z_3)][(z_1-z_3)+ t(z_2+z_3)]=0\, .
\]
It is straightforward to see that this family defines a lattice--isotopy from
$\A(\G)$ to $D_3$.
\end{proof}

\section{Proof of Theorems A, B and C}
\label{sect5}

We need one more ingredient: modular inequalities.

\subsection{}
\label{ss51}

These inequalities may be formulated for arbitrary connected CW--spaces 
of finite type, $M$, endowed with a $1$--marking, that is, a distinguished
$\Z$--basis of $H_1(M)$. The marking allows us to extend Definition \ref{def:ratl}
verbatim, to this more general context, as well as the definition of
$b_q(M, {\bf k}/d)$. 

Consider next the prime field $\k= \F_p$. In the presence of the marking,
we may speak about the element $\omega_{{\bf k}}\in H^1(M, \F_p)$, 
defined by taking the mod
$p$ reduction of ${\bf k}$. Hence, there is an associated Aomoto complex,
$(H^{\bullet}(M, \F_p), \mu_{{\bf k}})$, defined exactly as in
\eqref{eq:aok}, leading to the numbers $\beta_{qp}(M, {\bf k})$; see \eqref{eq:modres}.
When $M=M_{\A}$ is an arrangement complement, 
$\beta_{qp}(M_{\A}, {\bf 1})= \beta_{qp}(\A)$.

\begin{theorem}[\cite{PS}]
\label{thm:cohen}
Assume that the connected, finite type, $1$--marked CW--space $M$ 
has torsion--free integral homology.
Let $\rho$ be a rational local system on $M$, with denominator $d=p^s$, 
where $p$ is prime and $s\ge 1$. Then
\[
b_q(M, {\bf k}/d) \le \beta_{qp}(M, {\bf k})\, , \forall\, q\,.
\]
\end{theorem}

This extends a result from \cite{OC}, where $M$ is an arrangement complement, and $s=1$. 

\begin{corollary}
\label{cor:modb}
Let $\A$ be a central arrangement of $n$ hyperplanes, and $p$ be a prime such that
$d:= p^s$ divides $n$. If $\beta_{qp}(\A)=0$, for $q\le k$, then $b_{qd}(\A)=0$,
for $q\le k$.
\end{corollary}

\begin{proof}
By Theorem \ref{thm:cohen}, $b_q(\A, {\bf 1}/d)=0$, for $q\le k$. Hence,
$b_{qd}(\A)=0$, for $q\le k$, by \eqref{eq:twistedf} and induction.
\end{proof}

\subsection{}
\label{pfa}

{\bf Proof of Theorem \ref{thm:main}.}

{\em Part} \eqref{a1}. Use figure \ref{fig:graph34}  to infer
that the ${\bf m}_2$--list of $\A(\G)$ from \eqref{eq:mlist} must be 
contained in $\{3, 4\}$. Therefore, Theorem \ref{thm:negvanish} implies that
$b_d(\G)=0$, if $d\ne 2, 3, 4$. For $d=2$ or $4$, recall from Proposition \ref{p2}
that $\beta_2(\G)=0$, and use Corollary \ref{cor:modb} to obtain again the vanishing
of $b_d(\G)$, as asserted.

{\em Part} \eqref{a2}. Follows from Lemma \ref{pnot234} and Proposition \ref{p2}.

{\em Part \eqref{a3}}. By inspecting the graphs from 
Figures \ref{5graphs} and \ref{2moregraphs},
we deduce from Proposition \ref{p3}, in conjunction with  Corollary \ref{cor:modb}, 
that either $\G$ is not exceptional, and then 
$\beta_3(\G)=0$, hence
$b_3(\G)=\beta_3(\G)=0$ for $n \equiv 0$ (mod $3$), or
$\G$ is  exceptional, and then $n\equiv 0$ (mod $3$) and $\beta_3(\G)=1$. 
Therefore, the proof of Part \eqref{a3} is reduced, via Lemma \ref{liso},
to checking that $b_3(D_3)=b_3(D_4)=1$.

This in turn may be easily done by using the Deligne method, as follows.
Choose integers $a$, $b$ and $c$ such that $a+b+c=-1$. Next, set
$\alpha_{12}^{\pm}=\alpha_{34}^{\pm}=a$, $\alpha_{13}^{\pm}=\alpha_{24}^{\pm}=b$, and
$\alpha_{23}^{\pm}=\alpha_{14}^{\pm}=c$. View $\{ \alpha_{ij}^{\pm} \}_{1\le i<j\le v}$
as an element $\alpha\in A^1_{\Z}(D_v)$, for $v=3, 4$.
It is easy to verify that $\frac{{\bf 1}}{3}+\alpha \in A^1_{\C}(D_v)$ is
$1$--nonresonant, in the sense of Definition \ref{def=kgen}. Hence, Proposition
\ref{thm:partial gen} applies and gives that $b_3(D_v)= b_1(D_v, \frac{{\bf 1}}{3})=
\beta_1(D_v, \frac{{\bf 1}}{3}+\alpha)$. The Aomoto Betti number 
$\beta_1(D_v, \frac{{\bf 1}}{3}+\alpha)$ is then computed directly from the definition
\eqref{eq:baom}, by easy linear algebra, as explained in \cite[Lemma 3.3]{LY}.

{\em Part \eqref{a4}}. Let us inspect the equivariant decomposition of $H_1(F_{\G}, \Q)$ from
\eqref{eq:ciclo}. As recalled in the Introduction, the divisor $d=1$ contributes with
exponent $n-1$. No other divisors can contribute, excepting $d=3$, by Part \eqref{a1};
at the same time, $\beta_2(\G)=\beta_5(\G)=0$, by Part \eqref{a2}. Finally,
$b_3(\G)=\beta_3(\G)$, by Part \eqref{a3}.

\subsection{}
\label{pfb}

{\bf Proof of Theorem \ref{thm:maincompl}.} We know from \S \ref{pfa} above that
\[
H_{1}(F_{\G}, \mathbb{Q})=\big( \frac{\mathbb{Q}[t]}{t-1} \big)^{n-1}\oplus
\big( \frac{\mathbb{Q}[t]}{t^{2}+t+1} \big)^{\beta_3(\G)}\, .
\]
Theorem \ref{thm:maincompl} follows then from Proposition \ref{p3} and Lemma \ref{liso}.

\subsection{}
\label{pfc}

{\bf Proof of Theorem \ref{thm:sharp}.}
Theorem \ref{thm:cohen} predicts inequalities
\begin{equation}
\label{eq:modgr}
b_1(\A(\G), {\bf 1}/p^s) \le \beta_{p}(\G)\,, \quad {\rm for}\quad s\ge 1 \, ,
\end{equation}
at each prime $p$. We have to show that they all are actually equalities, if $s=1$.

In rank $\ge 3$, this follows from Theorem \ref{thm:main}\eqref{a2}--\eqref{a3}.

This is equally true for an arbitrary rank $2$ arrangement $\A$. Indeed, in 
this case one knows that
\begin{equation} 
\label{eq:v2}
b_1(\A, {\bf 1}/d)=
\begin{cases}
0, & {\rm if}\; d\nmid n ;\\
n-2 , & {\rm if}\; d\mid n ,
\end{cases}
\end{equation}
where $n=\abs{\A}$ and $d\ne 1$, see for instance \cite[Example 10.1]{Su}.
As an immediate consequence of Lemma \ref{rules}, we also have
\begin{equation} 
\label{eq:r2}
\beta_p(\A)=
\begin{cases}
0, & {\rm if}\; p\nmid n ;\\
n-2 , & {\rm if}\; p\mid n ,
\end{cases}
\end{equation}
for every prime $p$.
Our assertion follows then by comparing \eqref{eq:v2} (for $d=p$) and \eqref{eq:r2}.

The proof of Theorem \ref{thm:sharp} is thus completed.

\begin{remark}
\label{rk:square}
When $n=\abs{\A}$ is prime, equations \eqref{eq:v2} and \eqref{eq:r2} above
also show that the inequality \eqref{eq:modgr} may well be strict, if $s>1$.
\end{remark}

\bibliographystyle{amsplain}

\end{document}